\input amstex
\documentstyle{dowlenfpfreel}

\pagewidth{5.4 in}
\pageheight{7.63 in}
\hcorrection{0.438 in}
\vcorrection{0.5 in}

\nopagenumbers
\headline={\tenrm \hfil \folio}

\topmatter
\TITLE
 Every nonreflexive subspace of $L_1[0,1]$ \\ fails
 the fixed point property
 \endTITLE
\author P.N. Dowling and C.J. Lennard\endauthor
\affil Miami University and University of Pittsburgh\endaffil
\address {P.N. Dowling\newline
          Mathematics and Statistics Department\newline
          Miami University\newline
          Oxford, Ohio 45056\newline
          U.S.A.\newline
          \newline
          C.J. Lennard\newline
          Mathematics and Statistics Department\newline
          University of Pittsburgh\newline
          Pittsburgh, Pennsylvania 15260\newline
          U.S.A.}

%
%
\predefine\barov{\B}
\redefine\B{  {\bf B} }
\def\varsxn { (x_n)_{n=1}^{\infty} }
\def\oo{ \infty }
\def\N{ {\bf N}   }
\def\||{ \Vert }

\def\R{ {\bf R} }

\def\BSX{ ( X , \Vert \cdot \Vert_X ) }
\def\eps{ \varepsilon }

\def\tto{ \rightarrow }
\def\lam{ \lambda }

\def\and{ \text{ and }   }
\def\ntto{ \underset n \to \longrightarrow  }
\abstract{ The main result of this paper
 is that every non-reflexive subspace $Y$ of $L_{\ 1}[0,1]$ fails
 the fixed point property for closed, bounded, convex subsets $C$ of
 $Y$ and nonexpansive (or contractive) mappings on $C$. 
 Combined with a theorem of Maurey we get
 that for subspaces $Y$ of $L_{\ 1}[0,1]$, $Y$ is reflexive if and only if
 $Y$ has the fixed point property. For general Banach spaces the question
 as to whether reflexivity implies the fixed point property and the converse 
 question are both still open.}


\endtopmatter
\document

\subheading{Introduction}

 We first give a new example
 of a renorming of $\ell_1$ such that the positive part of the usual
 closed unit ball supports a fixed point free nonexpansive map. 
 The mapping is that of Lim [L], whose renorming of $\ell_1$
 provided the first such example. This example
 leads to a new fixed point free nonexpansive map
 with respect to the usual norm. A variation on this theme produces 
 a fixed point free contraction on a closed, bounded, convex set in 
 $\ell_1$ with its usual norm. Using this example we show that every
 nonreflexive subspace of $L_1[0,1]$ fails the fixed point property
 for nonexpansive mappings, proving the converse of a theorem of
 Maurey [M].
 In particular, the Hardy space $H^1$ on the unit circle must fail
 to have the fixed point property; which contrasts with Maurey's
 result in [M] that $H^1$ has the weak (and weak-star) fixed point property.
 We next provide an example of a reflexive subspace of $L_1[0,1]$
 that fails normal structure; which means that Maurey's 
 first above-mentioned theorem cannot be
 deduced from Kirk's theorem.
 In passing, we observe that a standard example shows that $\ell_1$ 
 with its usual norm fails
 the weak-star fixed point property and the weak-star Kadec-Klee property
 with respect to the predual
 $c$ of all convergent scalar sequences.

 We thank
 Brailey Sims, Mark Smith and Barry Turett for helpful discussions.
 The second author acknowledges the support of a University of Pittsburgh
 F.A.S. Reseach Grant.

\subheading{ Preliminaries}
Recall that $\ell_1$ is the Banach space of all scalar sequences
 $x = \varsxn$ for which $\|| x \||_1 := \sum_{n=1}^{\oo} | x_n | < \oo$.
 For each $n \in \N$, we let
 $e_n$ denote the $n$-th vector in the usual unit vector basis of $\ell_1$.
 $L_1[0,1]$ is the usual space of Lebesgue integrable functions (where
 almost everywhere equal functions are identified), with its usual norm.

 Let $\BSX$ be a Banach space. We say that $\BSX$ has the fixed point property
 if given any non-empty, closed, bounded and convex subset $C$ of $X$,
 every nonexpansive mapping $T : C \tto C$ has a fixed point. Here $T$
 is nonexpansive if $\|| Tx - Ty \||_X \le \|| x - y \||_X $ for all
 $x, y \in C$. Moreover, $T$ is a contraction if 
  $\|| Tx - Ty \||_X < \|| x - y \||_X $ for every $x, y \in C$ with
 $x \neq y$.
 Also, $\BSX$ has normal structure if every closed, bounded,
 convex subset $C$ of $X$ containing at least two points has radius less than
 its diameter. For such a $C$ and $x \in C$, rad$(x;C)$ equals
 $\sup_{y \in C} \|| x - y \||_X$, the radius of $C$ is
 rad$(C) := \inf_{x \in C} \text{rad}(x;C)$
 and the diameter of $C$ is
 diam$(C) := \sup_{x \in C} \text{rad}(x;C)$. Kirk [Ki] showed that in a
 reflexive Banach space, normal structure implies the fixed point property.
 If $X$ is a dual space, isometrically isomorphic to $Y^*$ for some Banach
 space $Y$, then $\BSX$ has the weak-star fixed point property (with
 respect to $Y$) if given a non-empty, weak-star compact, convex set $C$
 in $X$, every nonexpansive mapping on $C$ has a fixed point. Further,
 $X$ has the weak-star
 Kadec-Klee property (with respect to $Y$) if
 weak-star and norm convergence of sequences coincides in the unit sphere
 of $X$. The weak Kadec-Klee property on a Banach space $X$ is defined
 similarly.

\subheading{1. A new example}

 Let $C := \{ x \in \ell_1^+ : \sum_{n=1}^{\oo} x_n \le 1 \}$. $C$ is a 
 closed, bounded and convex set in $\ell_1$. Define $T : C \tto C$
 by
 $$ T : x \longmapsto \left(
 1 - \sum_{n=1}^{\oo} x_n , x_1 , x_2, \dots , x_k ,
 \dots \right)    \ \ \ .$$
 $T$ is clearly fixed point free and $\sum_{n=1}^{\oo} (Tx)_n = 1$ for all
 $x \in C$. 

 Lim [L] renormed $\ell_1$ with the equivalent norm $x \longmapsto
 \|| x^+ \||_1 \vee \|| x^- \||_1 $. $T$ is nonexpansive with respect
 to Lim's norm.

 Let us renorm $\ell_1$ in another way. Define $\|| \cdot \||$ on $\ell_1$
 by setting
 $$ \|| x \|| := \left| \sum_{n=1}^{\oo} x_n \right| + \sum_{n=1}^{\oo} 
 | x_n |  \ ,\forall \ x \in \ell_1 \ \ \ .$$
 Then for all $x$ and $y$ in $C$,
 $$ T x - T y = \left( \sum_{n=1}^{\oo} (y_n - x_n) , x_1 - y_1 , x_2 - y_2 ,
 \dots , x_k - y_k , \dots \right)  \ \ ; \ \text{and so}$$
 $$ \|| T x - T y \|| = \left| \sum_{n=1}^{\oo} (y_n - x_n) + \sum_{k=1}^{\oo}
 (x_k - y_k) \right| + \left| \sum_{n=1}^{\oo} (y_n - x_n) \right| 
 + \sum_{k=1}^{\oo} | x_k - y_k | = \|| x - y \|| \ \ .$$
  Thus, in summary, $T$ is a fixed point free isometry
 on the closed, bounded and convex set $C$, for the
 equivalent norm $\|| \cdot \||$ on $\ell_1$.

\subheading{2. A new nonexpansive map for the usual norm on $\ell_1$}

 Define the linear mapping $Q$ from $(\ell_1 , \|| \cdot \||)$ into
  $(\ell_1 , \|| \cdot \||)_1$ by 
 $$x \longmapsto \left(\sum_{n=1}^{\oo} x_n , x_1 , x_2 , \dots , x_k , \dots
 \right) \ \ .$$
 Note that $\|| Q x \||_1 = \|| x \|| $ for all $x \in \ell_1$ and the range 
 of $Q$ is $V := \{ y \in \ell_1 : y_1 = \sum_{n=2}^{\oo} y_n \} $. Moreover,
 $Q^{-1} : y \longmapsto (y_2 , y_3 , \dots , y_k , \dots )$.

 Now, let $K := Q(C) = \{ y \in \ell_1^+ : y_1 = \sum_{n=2}^{\oo} y_n \le 1 \}
 $ and define the mapping $S : K \tto K$ by $S := Q T Q^{-1}$. It is easy
 to check that for all $y \in K$,
 $$ S y = ( 1 , 1 - y_1 , y_2 , y_3 , \dots , y_k , \dots )   \ \ ; $$
 and that $S$ is a fixed point free isometry on the closed, bounded,
 convex set $K$ in $\ell_1$ with respect to the usual norm.
   
\subheading{3. A variation on this theme}

 Fix a sequence $(\eps_k)_{k=1}^{\oo}$ in $[0,1)$ with
 $\sum_{k=1}^{\oo} \eps_k < \oo$. Fix $\delta \in (0,1]$. Define
 $T : \ell_1 \tto \ell_1$ by
 $$ T : x \longmapsto \left( \delta \left(
 1 - \sum_{n=1}^{\oo} x_n \right) , (1 - \eps_1 ) x_1 , (1 - \eps_2 ) x_2,
 \dots , (1 - \eps_k ) x_k , \dots \right)    \ \ \ .$$
 We claim that $T$ is fixed point free. Indeed, suppose not. Then there is
 an $x \in \ell_1$ with $T x = x$. Thus
 $$ \delta \left( 1 - \sum_{n=1}^{\oo} x_n \right) = x_1  \ \ \ \and \ \ \ 
 (1 - \eps_k) x_k = x_{k+1} \ , \ \forall \ k \in \N  \ .$$
 Consequently, for each $k \in \N$,
 $ x_{k+1} = x_1 \prod_{j=1}^k (1 - \eps_j)$.
 There are two cases. If $x_1 = 0$ then $x_n = 0$ for all $n \in \N$.
 So, from the first equation above, $\delta = 0$; a contradiction. If 
 $x_1 \neq 0$, then $\prod_{j=1}^{\oo} (1 - \eps_j) = \lim_{k \tto
 \oo} x_{k+1} = 0$. This is precisely equivalent to $\sum_{j=1}^{\oo}
 \eps_j = \oo$; a contradiction.

 As before, we define 
 $C := \{ x \in \ell_1^+ : \sum_{n=1}^{\oo} x_n \le 1 \}$. Then $T$
 maps $C$ into $C$. Indeed, each $(T x)_n \ge 0$ and
 $ \sum_{n=1}^{\oo} (T x)_n = \delta \left( 1 - \sum_{n=1}^{\oo} x_n 
 \right) + \sum_{k=1}^{\oo} (1 - \eps_k) x_k \le 1$.

 We verify that $T$ is non-expansive with respect to the norm $\||
 \cdot \||$. Fix $x, y \in C$.
 $$ \align 
 T x - T y &= \bigg( \delta \sum_{n=1}^{\oo} (y_n - x_n) ,
 (1 - \eps_1) (x_1 - y_1) , (1 - \eps_2) (x_2 - y_2) ,
 \dots , \\  & (1 - \eps_k) (x_k - y_k) ,
 \dots \bigg)   \ ; \endalign$$
 and so
 $$\align \|| T x - T y \|| &=
 \left| \delta \sum_{n=1}^{\oo} (y_n - x_n) + \sum_{k=1}^{\oo}
 (1 - \eps_k) (x_k - y_k) \right|     \\
 & \ \ + \delta \left| \sum_{n=1}^{\oo} (y_n - x_n) \right| 
 + \sum_{k=1}^{\oo} (1 - \eps_k) | x_k - y_k |    \\
 &\le \left| \sum_{n=1}^{\oo} (1 - \delta) (x_n - y_n) \right|
     + \sum_{k=1}^{\oo} \eps_k |x_k - y_k|   \\
 & \ \ + \delta \left| \sum_{n=1}^{\oo} (y_n - x_n) \right| 
 + \sum_{k=1}^{\oo} (1 - \eps_k) | x_k - y_k | 
 = \|| x - y \||  \ . \endalign$$

 In summary, $C$ is a closed, bounded, convex set in $\ell_1$ and $T$
 is a fixed point free nonexpansive map on $C$ for the equivalent norm
 $\|| \cdot \||$ on $\ell_1$.
 Further, as in the previous section, $S := Q T Q^{-1}$ is a fixed
 point free nonexpansive map on the closed, bounded, convex set
 $K := Q(C)$, for the usual norm $\|| \cdot \||_1$ on $\ell_1$.

 Briefly consider what happens when $\delta < 1$ and
 $0 < \eps_k \le 1 - \delta$ for all $k \in \N$. 
 Then $T$ is not an isometry because
 $0 \in C$ and $\|| T x - T 0 \|| < \|| x \||$ for all $x \in C$.
 Similarly, $S$ is not an isometry.

\subheading{4. Fixed point free contractions in $\ell_1$}

 Fix $\delta \in (0,1)$ and 
 a sequence $(\eps_k)_{k=1}^{\oo}$ in $(0,1 - \delta)$ with
 $\sum_{k=1}^{\oo} \eps_k < \oo$. Again define
 $T : \ell_1 \tto \ell_1$ by
 $$ T : x \longmapsto \left( \delta \left(
 1 - \sum_{n=1}^{\oo} x_n \right) , (1 - \eps_1 ) x_1 , (1 - \eps_2 ) x_2,
 \dots , (1 - \eps_k ) x_k , \dots \right)    \ \ \ .$$
 From the previous section, we know that $T$ is fixed point free.

 Fix $\alpha \in \R$ and $q > 0$. Define 
 $$  D := \left\{ x \in \ell_1 : \sum_{k=1}^{\oo} (1 - \delta -
 \eps_k) x_k = \alpha  \ \ \and \ \ \sum_{j=1}^{\oo} |x_j| \le q
 \right\}  \ \ . $$
 Note that given $\alpha$ we may always choose $q$ so large that $D$
 is non-empty. Moreover, $D$ is a norm closed, convex and bounded set
 in $\ell_1$. 

 Let us specialize our discussion to the following case. It is
 sufficient for our purposes, although other choices of parameters
 will also suffice. Let $\delta := 1/2$ and $\eps_1 := 1/4$. Let
 $\alpha = 1 / 2$ and fix $q > 0$ so large that $D \ne \emptyset$. We
 introduce the auxilliary parameter $\eta$ and given our previous
 choices and what we will require below of $\eta$, it turns out that
 $\eta := 3/4$ is the best choice. Finally, we define the sequence
 $(\eps_k)_{k=2}^{\oo}$ by the equations
 $$ \left( \frac{1}{2} - \eps_{n+1} \right) (1 - \eps_n) - \frac{1}{8}
 = \eta \left( \frac{1}{2} - \eps_n \right)    \ , \ \forall \ n \in \N \
 . $$
 It is easy to inductively show that $0 < \eps_n < 1/2 = 1 - \delta$
 for each $n$, and that the above equations for $\eps_n, n \ge 2$ are
 equivalent to
 $ \eps_{n+1} = \dfrac{\eps_n}{ 4 (1 - \eps_n)}$, for all $n
 \in \N$.
 Thus $0 < \eps_{n+1} < \eps_n / 2$ for each $n$, and so
 $\sum_{n=1}^{\oo} \eps_n < \oo$. Consequently, $\delta$ and
 $(\eps_k)_{k=1}^{\oo}$ satisfy the hypotheses of the beginning of
 this section. In particular, the mapping $T : \ell_1 \tto \ell_1$
 is fixed point free. 

\proclaim{4.1 Proposition} $T$ maps $D$ into $D$.
\endproclaim

\demo{\bf Proof} Fix $x \in D$. Then
 $$ \sum_{k=1}^{\oo} \left( \frac{1}{2} - \eps_k \right) x_k =
 \frac{1}{2} \ \ \and \ \ \sum_{j=1}^{\oo} |x_j| \le q  \ \ . $$
 $$ Tx = \left( \frac{1}{2} \left(
 1 - \sum_{n=1}^{\oo} x_n \right) , (1 - \eps_1 ) x_1 , (1 - \eps_2 ) x_2,
 \dots , (1 - \eps_k ) x_k , \dots \right)    \ \ .$$
 Thus, we see that
 $$ \align \sum_{n=1}^{\oo} \left( \frac{1}{2} - \eps_n \right) (Tx)_n
 &= \left( \frac{1}{2} - \eps_1 \right) \frac{1}{2} \left( 1 - 
 \sum_{n=1}^{\oo} x_n \right)  \\
 & \ \ + \sum_{n=1}^{\oo} \left( \frac{1}{2} - \eps_{n+1} \right)
 (1 - \eps_n) x_n              \\
 &= \frac{1}{8} - \frac{1}{8} \sum_{n=1}^{\oo} x_n
  + \sum_{n=1}^{\oo} \left( \frac{1}{2} - \eps_{n+1} \right)        
 (1 - \eps_n) x_n              \\
 &= \frac{1}{8}
  + \sum_{n=1}^{\oo}\left[ \left( \frac{1}{2} - \eps_{n+1} \right)        
 (1 - \eps_n) - \frac{1}{8} \right] x_n           \\
 &= \frac{1}{8}
  + \sum_{n=1}^{\oo} \frac{3}{4} 
 \left( \frac{1}{2} - \eps_n \right) x_n
 = \frac{1}{8} + \frac{3}{4} \left( \frac{1}{2} \right) = \frac{1}{2}
 \ \ ; \endalign $$
 and so we have
 $  \sum_{n=1}^{\oo} \left( \frac{1}{2} - \eps_n \right) (Tx)_n 
 = \frac{1}{2}$. Further,
 $$ \align \sum_{n=1}^{\oo} \left| (Tx)_n \right| &=
   \left| \frac{1}{2} \left( 1 - \sum_{n=1}^{\oo} x_n \right) \right|
  + \sum_{n=1}^{\oo} (1 - \eps_n) |x_n|          \\
 &=  \left| \frac{1}{2} - \frac{1}{2} \sum_{n=1}^{\oo} x_n  \right|
  + \sum_{n=1}^{\oo} (1 - \eps_n) |x_n|          \\
 &=  \left|  \sum_{n=1}^{\oo} - \eps_n x_n  \right|
  + \sum_{n=1}^{\oo} (1 - \eps_n) |x_n|          \\
 &\le    \sum_{n=1}^{\oo} \eps_n |x_n|
  + \sum_{n=1}^{\oo} (1 - \eps_n) |x_n| 
 \le    \sum_{n=1}^{\oo} |x_n| \le q   \ \ . \endalign $$
 So,
 $  \sum_{n=1}^{\oo} \left| (Tx)_n \right| \le q$.
 Thus $Tx \in D$ for all $x \in D$.\qed
\enddemo 

\proclaim{4.2 Proposition} $T$ is a contraction on $D$ with respect to
 the norm $\|| \cdot \||$ on $\ell_1$; and moreover, for all $x, y \in
 D$,
 $$\|| Tx - Ty \|| = \frac{1}{2} \left| \sum_{n=1}^{\oo} (y_n - x_n)
 \right| + \sum_{n=1}^{\oo} (1 - \eps_n) |x_n - y_n|  \ \ . $$
\endproclaim

\demo{\bf Proof} Fix $x, y \in D$. Then 
 $$\align
 T x - T y &= \bigg( \frac{1}{2} \sum_{n=1}^{\oo} (y_n - x_n) ,
 (1 - \eps_1) (x_1 - y_1) , (1 - \eps_2) (x_2 - y_2) ,
 \dots , \\ & (1 - \eps_k) (x_k - y_k) ,
 \dots \bigg)   \ . \endalign$$
 Therefore,
 $$\align \|| T x - T y \|| &=
 \left| \frac{1}{2} \sum_{n=1}^{\oo} (y_n - x_n) + \sum_{k=1}^{\oo}
 (1 - \eps_k) (x_k - y_k) \right|     \\
 & \ \ + \frac{1}{2} \left| \sum_{n=1}^{\oo} (y_n - x_n) \right| 
 + \sum_{k=1}^{\oo} (1 - \eps_k) | x_k - y_k |    \\
 &= \left| \sum_{n=1}^{\oo} \left( \frac{1}{2} - \eps_n \right) x_n
   -   \sum_{n=1}^{\oo} \left( \frac{1}{2} - \eps_n \right) y_n \right|  \\
 & \ \ + \frac{1}{2} \left| \sum_{n=1}^{\oo} (y_n - x_n) \right| 
 + \sum_{k=1}^{\oo} (1 - \eps_k) | x_k - y_k |    \\
 &= \left| \frac{1}{2} - \frac{1}{2} \right|
 + \frac{1}{2} \left| \sum_{n=1}^{\oo} (y_n - x_n) \right| 
 + \sum_{k=1}^{\oo} (1 - \eps_k) | x_k - y_k |  \ \ . 
 \endalign$$ \qed
\enddemo

 The following key result is now easy to verify. We omit the simple
 calculations involved.

\proclaim{4.3 Theorem} Let $D$ and $T$ be as in Propositions 4.1 and
 4.2. Then $S := Q T Q^{-1}$ is a fixed point free contraction on the
 non-empty,
 closed, bounded, convex set $L:= Q(D)$, for the usual norm
 $\|| \cdot \||_1$ on $\ell_1$. Moreover,
 $$ L = \left\{ y \in \ell_1 : y_1 = \sum_{n=2}^{\oo} y_n  \ , \ \ 
 \sum_{k=1}^{\oo} \left( \frac{1}{2} - \eps_k \right) y_{k+1} = \frac{1}{2}
 \ \and \ \ \sum_{n=2}^{\oo} |y_n| \le q   \right\}  \ \ ; 
 \text{while} $$
 $$ Sy = \left(
 1, \frac{1}{2} (1 - y_1), (1 - \eps_1) y_2,  (1 - \eps_2) y_3,
 \dots, (1 - \eps_k) y_{k+1}, \dots \right)  \ \forall \ y \in L \ ; \ \and $$
$$ \|| Sy - Sz \||_1
 = \frac{1}{2} |z_1 - y_1| + \sum_{n=1}^{\oo} (1 - \eps_n) 
 |y_{n+1} - z_{n+1}| \ \forall \ y, z \in L  \ . $$
\endproclaim

\subheading{5. All non-reflexive subspaces of $L_1[0,1]$ fail the FPP}
 
\proclaim{5.1 Proposition} Let $(\eps_n)_{n=1}^{\oo}$ be precisely the
 scalar sequence from results 4.1, 4.2 and 4.3 above.
 Let $\BSX$ be a Banach space and $Y$ be a
 subspace of $X$ such that there exists a sequence
 $(v_n)_{n=1}^{\oo}$ in $Y$, a sequence $(u_n)_{n=1}^{\oo}$ in $X$ and 
 a sequence $(\gamma_n)_{n=1}^{\oo}$ in $(0, \oo)$ with the following
 properties.
 $$ \left\|| \sum_{n=1}^{N} t_n u_n \right\||_X = \sum_{n=1}^{N}
 |t_n|  \ \ , \forall \ \text{scalar sequences} \ t_1, \dots t_N \ ;
 \tag i $$
 (i.e.  $(u_n)_{n=1}^{\oo}$ is an isometric copy in $X$ of the usual
 unit vector basis of $\ell_1$).
 $$  \|| u_n - v_n \||_X < \gamma_n  \ \ , \ \forall \ n \in \N \
 ; \tag ii $$
 where each $\gamma_n$ is chosen so small that
 $(v_n)_{n=1}^{\oo}$ spans an isomorphic copy of $\ell_1$
 in $Y$.
 $$ \gamma_1, \gamma_2 < \frac{1}{2} \left( \frac{1}{2} \right)
 \ \ \and \ \ \gamma_{n+1}, \gamma_{n+2} < \frac{1}{2} \eps_n 
 \ \forall \ n \in \N  \ \ . \tag iii $$
 Then $(Y, \|| \cdot \||_X)$ fails the fixed point property for closed, 
 bounded, convex sets in $Y$ and nonexpansive (or contractive)
 mappings on them.
\endproclaim 

\demo{\bf Proof} Let $q$ be as in 4.1, 4.2 and 4.3 above. We will
 apply Theorem 4.3. Define
 $$ \align
 M := \bigg\{ & \sum_{n=1}^{\oo} y_n v_n :
 y \in \ell_1 \ , \ \  y_1 = \sum_{n=2}^{\oo} y_n  \ ,    \\
 & \ \sum_{k=1}^{\oo} \left( \frac{1}{2} - \eps_k \right) y_{k+1} = \frac{1}{2}
 \ \and \ \ \sum_{n=2}^{\oo} |y_n| \le q   \bigg\}  \ \ .
\endalign $$
$M$ is a non-empty, closed, bounded, convex set in $Y$, by hypothesis
 (ii) and Theorem 4.3. Define $R : M \tto M$ by
 $$ R \left( \sum_{n=1}^{\oo} y_n v_n \right) :=
     \sum_{n=1}^{\oo} (Sy)_n v_n  \ \ , $$
 where $S : L \tto L$ is as in Theorem 4.3. Clearly, $R$ maps $M$ into
 $M$ and $R$ is fixed point free.

 Moreover, for all $\sigma = \sum_{n} y_n v_n$ and
  $\tau = \sum_{n} z_n v_n$ in $M$, Theorem 4.3 gives us that
 $$\align 
 \|| R\sigma - R\tau \||_X &\le 
 \left\|| \sum_{n=1}^{\oo} ((Sy)_n - (Sz)_n) u_n \right\||_X +
 \left\|| \sum_{n=1}^{\oo} ((Sy)_n - (Sz)_n) (v_n - u_n) \right\||_X   \\
 &\le \sum_{n=1}^{\oo} |(Sy)_n - (Sz)_n| + \sum_{n=1}^{\oo}
 |(Sy)_n - (Sz)_n|  \gamma_n   \\
 &\le \sum_{n=1}^{\oo} |(Sy)_n - (Sz)_n| + \sum_{m=1}^{\oo}
 |y_m - z_m|  \gamma_{m+1}   \\
 &\le \frac{1}{2} |z_1 - y_1| +
 \sum_{n=1}^{\oo} (1 - \eps_n)|y_{n+1} - z_{n+1}|    \\
 & \ \ + |y_1 - z_1| \frac{1}{4}
  + \sum_{m=2}^{\oo} |y_m - z_m| \frac{1}{2} \eps_{m-1}   \\
 &= \frac{3}{4} |y_1 - z_1| + \sum_{n=1}^{\oo}
 \left(1 - \frac{\eps_n}{2} \right)  |y_{n+1} - z_{n+1}|   
 \ \ . \endalign $$
 Meanwhile, for the same $\sigma$ and $\tau$ in $M$,
 $$\align 
 \|| \sigma - \tau \||_X &\ge 
 \left\|| \sum_{n=1}^{\oo} (y_n - z_n) u_n \right\||_X -
 \left\|| \sum_{n=1}^{\oo} (y_n - z_n) (v_n - u_n) \right\||_X   \\
 &\ge \sum_{n=1}^{\oo} |y_n - z_n| - \sum_{n=1}^{\oo} |y_n - z_n| 
 \gamma_n   \\
 &\ge \sum_{n=1}^{\oo} |y_n - z_n| - |y_1 - z_1| \frac{1}{4}
  - \sum_{m=2}^{\oo} |y_m - z_m| \frac{1}{2} \eps_{m-1}   \\
 &= \frac{3}{4} |y_1 - z_1| + \sum_{n=1}^{\oo}
 \left(1 - \frac{\eps_n}{2} \right)  |y_{n+1} - z_{n+1}|   
 \ \ . \endalign $$
 So,
 $ \|| R\sigma - R\tau \||_X \le \|| \sigma - \tau \||_X $,
 for all $\sigma, \tau \in M $.
 Note that $R$ is, in fact, a contraction on $M$.\qed
\enddemo

 Proposition 5.1 leads directly to our main result.
  
\proclaim{5.2 Theorem} Every non-reflexive subspace $Y$ of $L_1[0,1]$,
 with its usual norm, fails the fixed point property for closed, 
 bounded, convex sets in $Y$ and nonexpansive (or contractive)
 mappings on them.
\endproclaim

\demo{\bf Proof} By the proof of the Kadec-Pe\l czynski theorem [K-P]
 (or see [D], Chapter VII),
 for $X := L_1[0,1]$ with its usual norm, sequences
 $(v_n)_{n=1}^{\oo}$ in $Y$, $(u_n)_{n=1}^{\oo}$ in $X$ and 
 $(\gamma_n)_{n=1}^{\oo}$ in $(0, \oo)$ exist that satisfy all the 
 hypotheses of Proposition 5.1 above.
\qed
\enddemo

Combining 5.2 with Maurey's theorem [M] allows us to state the fact below.

\proclaim{5.3 Theorem} Let $Y$ be a subspace of $L_1[0,1]$ with its usual
 norm. Then the following are equivalent.
 
(i) $Y$ is reflexive.

(ii) $Y$ has the fixed point property.

\endproclaim

\subheading{6. There is a reflexive subspace of $L_1[0,1]$ failing NS}

 Let us define the sequence of Rademacher functions $(r_n)_{n=1}^{\oo}$
 on the real line $\R$. $r_0$ is the characteristic function of $[0,1]$,
 $r_1 := \chi_{(0,1/2)} - \chi_{(1/2,1)}$
 and $r_n(t) := r_{n-1}(2t) + r_{n-1} (2 t - 1)$, for all $t \in \R$
 and for each $n \ge 2$. Henceforth we will restrict the domain of each
 $r_n$ to $[0,1]$.

 We define the subspace $X$ of $L_1[0,1]$ to be the closed linear span
 of the sequence  $(r_n)_{n=0}^{\oo}$; which is isomorphic to $\ell_2$
 by Khinchine's inequalities (see, for example, Diestel [D]).
 In particular, $X$ is reflexive. We will
 denote the usual norm on $L_1[0,1]$ by $\|| \cdot \||_1$.

 Let $C$ be the closed convex hull in $X$ of
 the sequence  $(x_n)_{n=1}^{\oo}$, where each $x_n := r_n + r_0$. 
 (The translation factor of $r_0$ is introduced to show that we may
 arrange for $C$ to consist of all non-negative functions, and to link
 in with a remark below concerning the weak Kadec-Klee property).
 The set $C$ is closed, bounded and convex.

 We now show that $x \longmapsto$ rad$(x;C)$ is constant on $C$; and so
 $X$ with the norm $\|| \cdot \||_1$ fails normal structure.
 Indeed, fix $x,y \in C$. We may suppose, without loss of generality that
 $x,y \in D:=$ the convex hull of $\varsxn$. So there exist non-negative
 real numbers
 $\alpha_1, \dots , \alpha_N$ and $\beta_1 , \dots , \beta_N$, with
 $\sum_{n=1}^N \alpha_n = 1$ and $\sum_{n=1}^N \beta_n = 1$, such that
 $x = \sum_{n=1}^N \alpha_n x_n$ and $y = \sum_{n=1}^N \beta_n x_n$.
 Then
 $$\|| x - y \||_1 =  \left\|| \sum_{n,m=1}^N \alpha_n \beta_m 
 ( r_n - r_m ) \right\||_1 \le \max_{1 \le n,m \le N} \|| r_n - r_m \||_1 =1 
 \ \ . $$
 On the other hand, $x_{N+1} \in C$ and it is straightforward to calculate that
 $$\align
 \|| x - x_{N+1} \||_1 &= \left\|| \sum_{n=1}^N \alpha_n r_n - r_{N+1}
 \right\||_1                                                             \\
 &= \frac{1}{2^{N+1}} \sum_{\eps \in \{-1,1 \}^{N+1}}
 | \eps_1 \alpha_1 + \dots + \eps_N \alpha_N + \eps_{N+1} |              \\
 &= \frac{1}{2^{N+1}}\  2 \sum_{\eps \in \{-1,1 \}^{N}}
 ( 1 + \eps_1 \alpha_1 + \dots + \eps_N \alpha_N )                       \\
 &= \frac{1}{2^{N+1}} \ 2 \sum_{\eps \in \{-1,1 \}^{N}}
 ( 1 ) = 1      \ \ .
\endalign $$

 We remark that the above Banach space $X$ is a rather natural example
 of a reflexive space that has the fixed point property but fails normal
 structure. Karlovitz [Ka2] showed that another equivalent
 renorming of $\ell_2$
 (due to R.C. James) has the fixed point property, while it fails normal
 structure.

 Also note that $\varsxn$
 provides an example which shows that $X$ fails the weak Kadec-Klee
 property.
 
\subheading{7. An old example revisited}

 For each $x \in \ell_1$,
 we let $x_0 = \sum_{n=1}^{\oo} x_n$.
 One of the preduals of $\ell_1$
 is $c$, the space of all convergent scalar sequences
 $\lam = (\lam_n)_{n=1}^{\oo}$ with the supremum norm $\|| \lam \||_{\oo}$.
 We let $\lam_0 = \lim_{n \tto \oo} \lam_n$.
 The duality is given by (see, for example, Banach [B]) :
 $$ \langle x, \lam \rangle = \sum_{n=1}^{\oo} x_{n+1} (\lam_n - \lam_0)
 + x_0 \lam_0  = \sum_{n=0}^{\oo} x_{n+1} \lam_n   \ \ \ .$$

 Consider the following well-known example.
 Let $W$ equal $\{ x \in \ell_1^+ : \sum_{n=1}^{\oo} x_n = 1 \}$
 and $F : W \tto W$ be given by $F x := ( 0, x_1 , x_2 , \dots , x_k , \dots )
 $. Then $F$ is a fixed point free isometry on the closed, bounded, convex set
 $W$ in $\ell_1$. Let $\lam_n := 1$ for all $n \in \N$. Then $\lam \in c$ and
 $ \langle x , \lam \rangle = \sum_{n=1}^{\oo} x_n$.
 So, $W$ is weak-star compact with respect to the predual $c$ of $\ell_1$.
 Thus $\ell_1$ with its usual norm
 fails the weak-star fixed point property with respect
 to its predual $c$.

 For comparison, we remark that using
 the map $T$ and set $C$ described in section 1, Lim [L] showed that
 $\ell_1$ with Lim's norm fails the weak-star
 fixed point property for nonexpansive maps with respect to the
 isometric predual $c_0$ (suitably renormed). On the other hand,
 Karlovitz [Ka1] demonstrated
 that $\ell_1$ with its usual norm has the weak-star fixed point property
 with respect to $c_0$. Moreover, Soardi [So] showed that every Banach
 space $Y$ that is isomorphic to $\ell_1$, with Banach-Mazur distance from
 $\ell_1$ less than 2, has the weak-star fixed point property.

 In a similar manner to that in section 6, one can show that $W$
 is diametral, i.e. $x \longmapsto$ rad$(x;C)$ is constant on $C$.
 Moreover,
 the extreme points of $W$ provide us with an example showing that $\ell_1$
 fails the weak-star Kadec-Klee property with respect to the predual $c$.
 Let $x^{(k)} := e_{k+1}$ for each integer $k \ge 0$. Fix $\lam \in c$.
 For all $n \in \N$,
 $$ \langle x^{(n)}, \lam \rangle = \sum_{m=0}^{\oo} x^{(n)}_{m+1} \lam_m
   = \lam_n \ntto \lam_0 = \langle x^{(0)}, \lam \rangle \ \ . $$
 So $x^{(n)} \ntto x^{(0)}$ weak-star with respect to $c$, each
 $\|| x^{(k)} \||_1 = 1$, yet $\|| x^{(n)} - x^{(0)} \||_1 = 2$ for all 
 $n \ge 1$.

 We contrast the above with the fact that $\ell_1$ (with its usual norm)
 has the weak-star uniform
 Kadec-Klee property with respect to the predual $c_0$ (see, for example,
 [D-S] and [Si]); and therefore has both the weak-star fixed point property
 and the weak-star Kadec-Klee property with respect to $c_0$.

 Finally, we remark that Michael Smyth [Sm] has recently shown, by a variation
 of the techniques in section 3 and 4 above, that $\ell_1$ fails
 the weak-star fixed point property with respect to its predual $c$
 with a contractive map. By modifying our set $D$ in section 4, he is also
 able to slightly simplify our mapping $T$.
 
\newpage 

\heading References \endheading
\bigskip

[B] S. Banach, Th\'eorie des Op\'erations Lin\'eaires, Chelsea Publishing
 Co., New York 1978.

[D] J. Diestel, Sequences and Series in Banach Spaces, Springer-Verlag
 New York Inc. 1984.

[D-S] D. van Dulst and B. Sims, Fixed points of nonexpansive mappings
 and Chebyshev centers in Banach spaces with norms of type (KK), in:
 Banach Space Theory and its Applications, Proc. Bucharest 1981, Lecture
 Notes in Math. 991, Springer-Verlag 1983, 35-43.

[K-P] M.I. Kadec and A. Pe\l czynski, Bases, lacunary sequences
 and complemented subspaces in the spaces
 $L_p$, Studia Math. 21 (1962), 161-176.

[Ka1] L.A. Karlovitz, On nonexpansive mappings, Proc. Amer. Math. Soc.
 55(2) (1976), 321-325.

[Ka2] L.A. Karlovitz, Existence of fixed points of nonexpansive mappings
 in a space without normal structure, Pacific J. Math.
 66(1) (1976), 153-159.

[Ki] W.A. Kirk, A fixed point theorem for mappings which do not increase
 distances, Amer. Math. Monthly 72 (1965), 1004-1006.

[L] T.C. Lim, Asymptotic centers and nonexpansive mappings in conjugate
 Banach spaces, Pacific J. Math. 90(1) (1980), 135-143.

[M] B. Maurey, Points fixes des contractions de certains faiblement
 compacts de $L^1$, Seminaire d'Analyse Fonctionelle, Expos\'e no. VIII,
 \'Ecole Polytechnique,
 Centre de Math\'ematiques (1980-1981).
 
[Si] B. Sims, The existence question for fixed points of nonexpansive maps,
 Lecture Notes, Kent State Univ. 1986.

[Sm] M. Smyth, preprint, University of Newcastle, Australia.

[So] P.M. Soardi, Schauder bases and fixed points of nonexpansive mappings,
 Pacific J. Math. 101(1) (1982), 193-198.

\enddocument